\documentclass[12pt,a4paper,twoside,final,notitlepage, leqno]{article}
\usepackage[english]{babel}
\usepackage[T1]{fontenc}  
\usepackage{epsfig, graphicx, amssymb}
\usepackage{amsmath,amsthm,epsfig,amsfonts,bbm}  
\usepackage[bottom]{footmisc}
 
\def \smb {{\scriptstyle \bullet }}
\newcommand{\monitem}{ \smallskip \noindent $\bullet$ \quad  } 
\newcommand{\moneq}{\vspace*{-6pt} \begin{equation} \displaystyle } 
\newcommand{\moneqstar}{\vspace*{-6pt} \begin{equation*} \displaystyle } 
\newcommand{\monendstar}{\vspace*{-6pt} \end{equation*}   }
\newcommand{\monend}{\vspace*{-6pt} \end{equation}   }

\def\section*#1{}

\usepackage{fancyhdr}
\renewcommand{\headrulewidth}{0pt}

\begin{document} 

\fancypagestyle{plain}{ \fancyfoot{} \renewcommand{\footrulewidth}{0pt}}
\fancypagestyle{plain}{ \fancyhead{} \renewcommand{\headrulewidth}{0pt}}

~


 \bigskip   \bigskip   \bigskip  
 
\centerline {\bf \LARGE  Pontryagin calculus in Riemannian geometry }

\bigskip \bigskip \bigskip \bigskip 

\centerline { \large Fran\c cois Dubois$^{a,b}$, Danielle Fortun\'e$^{c}$, Juan Antonio Rojas Quintero$^{d}$}
\smallskip 
\centerline { \large and Claude Vall\'ee}    

\bigskip  \bigskip \bigskip 

\centerline { \it  \small   
$^a$    Conservatoire National des Arts et M\'etiers,} 

\centerline { \it  \small 
Laboratoire de M\'ecanique des Structures 
et des Syst\`emes Coupl\'es, F-75003,  Paris, France}

\centerline { \it  \small 
$^b$  Department of Mathematics, University  Paris-Sud, } 

\centerline { \it  \small 
 B\^at. 425, F-91405 Orsay Cedex, France.}

\centerline { \it  \small 
$^c$  21, rue du Hameau du Cherpe, 86280 Saint Beno\^{\i}t, France } 

\centerline { \it  \small 
$^d$  School of Mechanics and Engineering, \\ 
Southwest Jiaotong University, Chengdu, China }

\bigskip   
\centerline { \it  \small francois.dubois@math.u-psud.fr, danielle.fortune123@orange.fr, jrojasqu@yahoo.com }

\bigskip   

\bigskip \bigskip 

\centerline { {\rm  03 June 2015 } 
\footnote {\rm  \small $\,$ 
Contribution published in {\it 
Geometric Science of Information - Second International Conference, GSI 2015}, 
Editors: F. Nielsen, F. Barbaresco, 
Series: Lecture Notes in Computer Science, volume~9389, 
Springer,  pages~541-549, november 2015. }}

\bigskip  \bigskip \bigskip 
\noindent  {\bf Abstract. } \qquad 
In this contribution, we study systems 
 with a finite number 
of degrees of freedom as in robotics.  
A key idea is to consider the mass tensor associated
to the kinetic energy as a metric in a Riemannian configuration space. 
We apply Pontryagin's framework to derive 
an optimal evolution of the control forces and torques 
applied to the mechanical system. 
This equation under covariant form uses explicitly the Riemann curvature tensor.

\noindent 
This contribution is dedicated to the memory of Claude Vall\'ee (1945-2014).

 $ $ \\  [2mm]

\bigskip \bigskip  \newpage \noindent {\bf \large   Introduction}  
%
  \fancyfoot[C]{\oldstylenums{\thepage}}
\fancyhead[EC]{\sc{F. Dubois, D. Fortun\'e, J.A. Rojas Quintero and C. Vall\'ee}} 
\fancyhead[OC]{\sc{Pontryagin  calculus in Riemannian geometry}} 
\fancyfoot[C]{\oldstylenums{\thepage}}
%
  
\fancypagestyle{plain}{ \fancyfoot{} \renewcommand{\footrulewidth}{0pt}}
\fancypagestyle{plain}{ \fancyhead{} \renewcommand{\headrulewidth}{0pt}} 

\bibliographystyle{alpha}


\noindent 
As part of studies based on the calculus of variations, 
the choice of a Lagrangian or a Hamiltonian is essential.
When we study the dynamics of articulated systems, the choice 
of the Lagrangian is directly linked to the conservation of energy. 
The Euler-Lagrange methodology is applied to the kinetic and potential energies. 
This establishes a system of second order ordinary differential equations for the  movement.
These equations are identical to those deduced from the fundamental principle of dynamics. 
The choice of configuration parameters does not affect  the energy value.
 Because the kinetic energy is a positive definite quadratic form 
with respect to the configuration parameters derivatives, its coefficients 
are ideal candidates to define and create a Riemannian metric structure 
on the configuration space. 
The Euler Lagrange equations have a contravariant tensorial nature 
and highlight the covariant derivatives with respect to time with 
the introduction of the Christoffel symbols. 

\monitem 
For the control of articulated robot choosing a Hamiltonian 
and a cost function are delicate. Here the presence of the 
Riemann structure is sound. It enables a cost function invariant 
when coordinates change. The application of the Pontryagin method 
from the optimal Hamiltonian leads to a system of second order 
differential equations for the control  variables. Its tensorial 
nature is covariant and the Riemann-Christoffel curvature tensor  
is naturally revealed. In this development, the adjoint variables 
are directly interpreted and have a physical sense.

\bigskip   \noindent {\bf \large 1) \quad Pontryagin  framework for differential equations} 

\noindent 
We study a  dynamical system, 
where the state vector   $\, y(t\,;\, \lambda(\smb)) \,$ 
is a function of time. This system is controlled by a set of variables  $\,\lambda(t) \,$ 
and satisfies a first order ordinary differential equation: 
\moneq \label{edo} 
{{{\rm d}y}\over{{\rm d}t}} \,\,=\,\, f(y(t) ,\, \lambda(t) ,\, t) \,.
\monend 
We suppose also given an initial condition:  
\moneq \label{condini} 
y(0\,;\,\lambda(\smb))  =  x \,. 
\monend 
We search an optimal solution associated to the optimal control
$\, t \longmapsto \lambda(t)  \,$  in order to minimize the 
following cost function $\, J \,$: 
\moneq  \label{fonction-cout} 
J( \lambda(\smb)) \, \equiv \, 
\int_{0}^{T} \, g \big( y(t),\, \lambda(t),\,t \big) \, {\rm d}t \,\,, 
\monend
where  $\, g(\smb ,\, \smb ,\, \smb) \, $ is a given real valued function.

\monitem 
Pontryagin's main idea \cite{PBGM62} is to consider the differential equation  
(\ref{edo}) 
as a {\it constraint} satisfied by the variable  $\, y $. 
Then he  introduces a Lagrange multiplier $\, p \,$
associated to the constraint  (\ref{edo}). This new variable, due to the continuous nature of the 
constraint  (\ref{edo}), is a covariant  vector function of  time:  $\, p = p(t) .\,$ 
A global Lagrangian functional can be considered:
\moneqstar 
{\cal L}(y ,\,\lambda ,\,p) \,\, \equiv \,\,
\int_{0}^{T} \, g(y,\,\lambda,\,t) \, {\rm d}t \,\, + \,\, \int_{0}^{T} \, p \, \Bigl(
{{{\rm d}y}\over{{\rm d}t}} - f(y ,\, \lambda ,\, t) \Bigr) \, {\rm d}t \,.\,
\monendstar

 \bigskip \monitem   {\bf Proposition 1. \quad  Adjoint equations. }

\noindent 
If the Lagrange multiplier   $\, p(t) \,$ satisfies the so-called {\it adjoint equations}, 
\moneq \label{eqs-adjointes} 
{{{\rm d}p}\over{{\rm d}t}} +  p \, {{\partial f}\over{\partial y}}   
- {{\partial g}\over{\partial y}}  \,\,=\,\, 0 \,
\monend 
and the so-called  {\it final condition},  
\moneq \label{cond-finale} 
p(T) \,\,=\,\, 0  \,,\,
\monend  
then the variation of the cost function 
for a given variation $\, \delta \lambda \,$ of the paramater 
is given by the simple relation 
\moneqstar 
\delta J \,=\,  \int_{0}^{T} \, \Bigl[ \,  {{\partial g}\over{\partial  \lambda}} - p \, {{\partial
f}\over{\partial \lambda}}  \, \Bigl] \,  \delta \lambda (t)   \, \,  {\rm d}t \,. 
\monendstar
At the optimum this variation is identically null and we find the so-called
Pontryagin optimality condition:
\moneq \label{cond-optimale}   
{{\partial g}\over{\partial  \lambda}} - p \, {{\partial f}\over{\partial \lambda}}  
\,=\, 0 \, . 
\monend  
%

%
\smallskip \noindent  {\bf Proof of Proposition 1.}

\noindent 
We write in a general way the variation of the Lagrangian  $\, {\cal L}(y ,\, \lambda ,\,p) \,$ 
in a variation  $\, \delta y ,\,$ $\, \delta \lambda \,$ and  $\, \delta p \,$
of the variables $\, y ,\, \lambda ,\, $ and  $\, p \,$  respectively.
We use classical  calculus rules as 
\moneqstar   
\delta \Big( \int_0^T \, g \, {\rm d}t \Big) \,=\, \int_0^T \, \delta g \, {\rm d}t \,, \quad 
\delta \big( {{{\rm d}y}\over{{\rm d}t}} \big) \,=\,
 {{{\rm d}}\over{{\rm d}t}} \big( \delta y \big) \,, 
\monendstar 
and we integrate by parts. We get 

\smallskip \noindent $  \displaystyle
\delta {\cal L} \,\,=\,\, 
\int_{0}^{T} \, \Bigl[ \,  {{\partial g}\over{\partial y}} \, \delta y + 
 {{\partial g}\over{\partial  \lambda }} \, \delta \lambda  \, \Bigl] \,  {\rm d}t \,\,
 + \,\, \int_{0}^{T} \,  p
\, \Bigl({{{\rm d} \delta y}\over{{\rm d}t}} - {{\partial f}\over{\partial y}} \delta y - 
 {{\partial f}\over{\partial  \lambda }} \delta \lambda \Bigr) \, {\rm d}t $ 

 \noindent $  \displaystyle \hfill 
+ \,\, \int_{0}^{T} \, \delta p \, \Bigl(
{{{\rm d}y}\over{{\rm d}t}} - f(y ,\, \lambda ,\, t) \Bigr) \, {\rm d}t \,\,$

\smallskip \noindent $  \displaystyle 
=\,\, 
\int_{0}^{T} \, \Bigl[ \,  {{\partial g}\over{\partial y}} - p \, {{\partial
f}\over{\partial y}}  \, \Bigl] \,  \delta y  \,  {\rm d}t \,  \,+\, 
\int_{0}^{T} \, \Bigl[ \,  {{\partial g}\over{\partial  \lambda }} - p \, {{\partial
f}\over{\partial  \lambda }}  \, \Bigl] \,  \delta \lambda  \,  {\rm d}t \, 
+ \,\, \Bigl[ \, p \, \delta y \, \Bigr]_{0}^{T} \, - \, \int_{0}^{T} \,
{{{\rm d}p}\over{{\rm d}t}} \,\, \delta y \, {\rm d}t \, $

\smallskip \noindent $  \displaystyle
=\,\, 
  p(T) \, \delta y(T) \, - \,  \int_{0}^{T} \,  \Bigl[ \, {{{\rm d}p}\over{{\rm d}t}} 
+ p \, {{\partial f}\over{\partial y}}  -  {{\partial g}\over{\partial y}}  \, \Bigl] \, 
\delta y  \, {\rm d}t \, 
+ \,\, \int_{0}^{T} \, \Bigl[ \,  {{\partial g}\over{\partial  \lambda }} - p \, {{\partial
f}\over{\partial  \lambda }}  \, \Bigl] \,  \delta \lambda  \, \,  {\rm d}t \,    $ 

\smallskip \noindent
because  $\, \delta y(0) = 0 \,$ taking fixed the initial condition  (\ref{condini}). 
By canceling the first two terms of the right hand side of the previous relation, we find
the adjoint equation (\ref{eqs-adjointes}) giving the evolution of the Lagrange multiplier
and the associated final condition (\ref{cond-finale}). The third term allows to calculate the
change in the functional $ \, J (\smb) \, $ for a  variation $ \, \delta  \lambda $ of
control.~$\square$

\bigskip    
\noindent {\bf \large 2) \quad Pontryagin  hamiltonian } 

\noindent 
We introduce the Hamiltonian 
\moneq \label{hamiltonien} 
 {\cal H} (p, \, y, \, \lambda)  \, \equiv \, p \, f - g 
\monend
and the optimal Hamiltonian 
\moneqstar 
 H (p, \, y) \equiv   {\cal H}  (p, \, y, \, \lambda^*) 
\monendstar 
for  $ \, \lambda (t) =  \lambda^* (t)\,$ equal to the optimal value associated to the optimal 
condition (\ref{cond-optimale}). 

\bigskip   \monitem    {\bf Proposition 2. \quad  Symplectic form of the dynamic equations. }

\noindent 
With the notations introduced proviously, the ``forward'' differential equation
(\ref{edo}) and the ``backward'' adjoint differential equation (\ref{eqs-adjointes}) take the 
following symplectic form:
\moneq \label{symplectic} 
{{{\rm d}y}\over{{\rm d}t}} \,= \, {{\partial H}\over{\partial p}}  \,,
\qquad 
{{{\rm d}p}\over{{\rm d}t}}\,= \, - {{\partial H}\over{\partial y}}  \,. 
\monend
%

%
\smallskip \noindent  {\bf Proof of Proposition 2.}

\noindent 
Since $ \, {{\partial {\cal H}}\over{\partial \lambda}} = 0 \, $ at 
the optimum, we have 
$ \,\, 
{{\partial H}\over{\partial p}} \,=\,  {{\partial {\cal H}}\over{\partial p}} 
\,=\, f \,$ and the first relation of
(\ref{symplectic}) is proven. On the other hand, 
$ \,\,  {{\partial H}\over{\partial y}} \,=\,  {{\partial {\cal H}}\over{\partial y}} $ 
$ = \, p \,  {{\partial f}\over{\partial y}} -  {{\partial g}\over{\partial y}} $ 
and the property is established. \hfill $\square $

\bigskip   \noindent {\bf \large 3) \quad Riemanian metric } 

\noindent 
We consider now a dynamical system parameterized by 
a finite number of functions $ \, q^j(t) \,  $
like a poly-articulated system for  robotics applications, 
as developed previously in \cite{LV95,RQ13,RVGSA13,VRFG13}. 
The set of all states $ \,  q \, \equiv \, \{ q^j \}  \, $ is denoted by $\, Q$. 
The kinetic energy $ \, K \, $
is a positive definite quadratic form of the time derivatives $ \, {\dot q}^j \,$
for each state $ \, q \in Q $.   The coefficients
of this quadratic form define a so-called {\it  mass tensor} $ \, M(q)  $.
The mass tensor 
is composed by  {\it a priori } a nonlinear regular function 
of the state $ \, q \in Q $.  We have 
\moneq \label{energie-cinetic} 
K(q, {\dot q}) \, \equiv \, {1\over2} \, \sum_{k \, \ell}  M_{k  \ell} (q) \, {\dot q}^k \, {\dot q}^\ell \, . 
\monend
The mass tensor $\, M(q) \,$ in  (\ref{energie-cinetic}) is symmetric 
and positive definite   for each state~$ q$. 
It contains the mechanical caracteristics of mass, inertia of the articulated system. 
With  Lazrak and  Vall\'ee \cite{LV95} and Siebert \cite{Si12}, we 
consider the Riemannian metric $ \, g \,$  defined by the mass  tensor  $M$.
We set:
\moneqstar 
g_{k  \ell}  (q) \, \equiv \,  M_{k  \ell}  (q) \, . 
\monendstar 

\monitem 
With this framework, the space of states $ \, Q \,$ 
has now a structure of {\it Riemannian manifold}.
Therefore all  classical geometrical tools of Riemannian geometry can be used 
(see {\it e.g.} the book \cite{LR75}): 

covariant space derivation $\, \partial_j \equiv {{\partial}\over{\partial q^j}}$ 

contravariant space derivation $ \,\, \partial^j $: $ \, \, < \partial^j \,,\,   \partial_k > 
\,=\, \delta^j_k $ 

component $ \, j, \ell \,$ of the inverse mass tensor $ \, M^{-1} $:  $ \, M^{j \ell} $,  
$ \,\, M_{i j} \,  M^{j \ell}  \,=\, \delta^\ell_i $ 

connection $\, \Gamma^j_{i k} \,=\,  {1\over2} \, M^{j \ell} \, \big( \partial_i M_{\ell k} 
+ \partial_k M_{\ell i} - \partial_\ell M_{i  k} \big) $,  
$ \, \, \Gamma^j_{k i} \,=\,\Gamma^j_{i k} $, 

\qquad \qquad $ \, {\rm d} \, \partial_j \,=\, \Gamma^\ell_{j k} \, {\rm d}q^k \, \partial_\ell $, 
$ \, {\rm d} \, \partial^j \,=\, - \Gamma^j_{k \ell} \, {\rm d}q^k \, \partial^\ell $, 

relations between covariant components $\, \varphi_j \,$
and contravariant components~$\, \varphi^k \,$ 

\qquad \qquad
of a vector field: $ \,\,  \varphi_j \,= \, M_{j k} \,  \varphi^k $, 
$ \,  \varphi^k \,= \, M^{k j} \,  \varphi_j $ 

covariant  derivation of a vector field $ \, \varphi \equiv \varphi^j \, \partial_j $:
$ \, {\rm d} \varphi  \, = \,  
\big( \partial_\ell \varphi^j +  \Gamma^j_{\ell k} \,  \varphi^k  \big) \, {\rm d}q^\ell \, \partial_j $

covariant  derivation of a covector field $ \, \varphi \equiv \varphi_\ell \, \partial^\ell $:

\qquad \qquad
$ \, {\rm d} \varphi  \, = \,  
\big( \partial_k \varphi_\ell -  \Gamma^j_{k \ell} \,  \varphi_j  \big) \, {\rm d}q^k \, \partial^\ell $

Ricci identities: 
%
%
\moneq \label{ricci} 
 \left\{ \begin{array}{rl}  \displaystyle \,\,  
  \partial_j M_{k \ell}  & = \, \, 
\Gamma^p_{j k}  \,  M_{\ell p} +  \Gamma^p_{j \ell}  \,  M_{k  p}
\,, \\    \displaystyle  \,\,  
\partial_j M^{k \ell} & = 
\, - \Gamma^k_{j p} \, M^{p \ell}  - \Gamma^\ell_{j p} \, M^{p k} 
\end{array}  \right. \monend

gradient of a scalar field: $ \, {\rm d} V \,=\, \partial_\ell V \, {\rm d}q^\ell  \,=\,
< \nabla V \,,\, {\rm d}q^j \, \partial_j > \, $ and  

\qquad \qquad
$ \, \nabla V \,=\,  \partial_\ell V \, \partial^\ell  $ 

gradient of a  covector field $ \, \varphi \,=\,  \varphi_\ell \, \partial^\ell $:
$ \, {\rm d} \varphi \,\equiv \, < \nabla \varphi \,,\,  {\rm d}q^j \, \partial_j > \, $
and 

\qquad \qquad
$ \, \nabla \varphi \,=\, 
\big( \partial_k \varphi_\ell -  \Gamma^j_{k \ell} \,  \varphi_j  \big) \, \partial^k \,  \partial^\ell $ 

second order gradient of a scalar field $V$: $ \, \nabla^2 V \,=\, \nabla ( \nabla V) \, $ and  
\moneq \label{second-gradient}  
\nabla^2 V \,=\, 
\big( \partial_k \partial_\ell V -  \Gamma^j_{k \ell} \,  \partial_j V \big) \, \partial^k \,  \partial^\ell 
\monend 

components $ \,  R^j_{i k \ell} \, $  of the Riemann tensor: 
\moneq \label{tenseur-riemann} 
  R^j_{i k \ell} \, \equiv \, \partial_\ell \Gamma^j_{i k} - \partial_k \Gamma^j_{i \ell}
+ \Gamma^p_{i k} \, \Gamma^j_{p \ell} - \Gamma^p_{i \ell} \, \Gamma^j_{p k}
\monend 

anti-symmetry of the Riemann tensor: 
$ \, R^j_{i k \ell} \, = \, - R^j_{i \ell k} $.

\bigskip   \monitem   {\bf Proposition 3. \quad  
Riemanian form of the Euler-Lagrange equations. }

\noindent 
With the previous framework, in the presence of an external potential  $ \, V = V(q)  $, 
the Lagrangian $\, L(q, \, \dot q) \, = \, K(q, {\dot q}) - V(q) \, $
allows to write the equations of motion in the classical   Euler-Lagrange form:
\moneq \label{euler-lagrange} 
{{\rm d}\over{{\rm d}t}} \Big( {{\partial L}\over{\partial {\dot q}^i}} \Big) 
\,= \,  {{\partial L}\over{\partial q^i}} 
\monend
These equations take also the Riemannian form: 
\moneq  \label{riemann-motion} 
M_{k \ell} \, \big( {\ddot q}^\ell + \Gamma^\ell_{i j} \, {\dot q}^i \,{\dot q}^j \big) \,  
+ \, \partial_k V \,= \, 0 \, . 
\monend

%
\smallskip \noindent  {\bf Proof of Proposition 3.}

\noindent 
The proof is presented in the references \cite{RVGSA13} and \cite{VRFG13}. 
We detail it to be complete. We have, due to (\ref{energie-cinetic}), 

\moneqstar 
{{\partial K}\over{\partial {\dot q}^k}} \,=\, M_{k \ell} \, {\dot q}^\ell  \, . 
\monendstar 
We have also the following calculus:

\smallskip \noindent $ \displaystyle 
{{\rm d}\over{{\rm d}t}} \Big( {{\partial L}\over{\partial {\dot q}^k}} \Big) 
-  {{\partial L}\over{\partial q^k}}  \, = \, 
{{\rm d}\over{{\rm d}t}}  \Big( {{\partial K}\over{\partial {\dot q}^k}} \Big) - 
\partial_k \big( K - V(q) \big) $

\smallskip \noindent $ \displaystyle \quad = \, 
{{\rm d}\over{{\rm d}t}} \big(M_{k \ell} \, {\dot q}^\ell \big) 
- \partial_k \Big( {1\over2} \,  M_{i j} \, {\dot q}^i \, {\dot q}^j \Big) 
+  \partial_k V $ 

\smallskip \noindent $ \displaystyle \quad = \, 
\big( \partial_j M_{k \ell} \big) \,  {\dot q}^j \, \, {\dot q}^\ell 
+  M_{k \ell}  \, {\ddot q}^\ell  
-  {1\over2}  \, \big( \partial_k M_{i j} \big) \,  {\dot q}^i \, {\dot q}^j
+  \partial_k V $ 

\smallskip \noindent $ \displaystyle \quad = \, 
\big( M_{k s} \, \Gamma^s_{j \ell} + M_{\ell s} \, \Gamma^s_{j k} \big) \,  {\dot q}^j \, \, {\dot q}^\ell 
 -  {1\over2}  \, \big(  M_{i s} \, \Gamma^s_{k j} + M_{j s} \, \Gamma^s_{k i} \big)\,  {\dot q}^i \, {\dot q}^j
+  M_{k \ell} \,  {\ddot q}^\ell  
+  \partial_k V $ 

\hfill due to the first Ricci identity  (\ref{ricci})

\noindent $ \displaystyle \quad = \, 
\Big(  M_{k s} \, \Gamma^s_{j i} +  M_{i s} \, \Gamma^s_{j k}  \, 
 -  {1\over2}  \,  M_{i s} \, \Gamma^s_{k j}  -  {1\over2}  \, M_{j s} \, \Gamma^s_{k i} \Big) 
\,  {\dot q}^i \, {\dot q}^j  
+  M_{k \ell} \,  {\ddot q}^\ell +  \partial_k V $ 

\smallskip  \noindent $ \displaystyle \quad = \, 
M_{k \ell} \, \Big( {\ddot q}^\ell +  \Gamma^\ell_{i j} \,  {\dot q}^i \, \, {\dot q}^j \Big) \,  +  \partial_k V $ 

\smallskip  \noindent 
and the relation (\ref{riemann-motion}) is established.  \hfill $\square$ 

\smallskip \monitem
When a mechanical forcing control $ \, u \, $ is present
(forces and torques typically),  the equations of motion can be formulated 
as follows:
\moneq \label{control-evolution} 
 {\ddot q}^j +  \Gamma^j_{k \ell} \,  {\dot q}^k \,{\dot q}^\ell +
 M^{j \ell} \, \partial_\ell V \, =  \, u^j \, . 
\monend
We observe that with this form (\ref{control-evolution}) of the equations
of motion, the contravariant components of the control $u$ have to be considered 
in the right hand side of the dynamical equations. 

\bigskip   \noindent {\bf \large 4) \quad Optimal dynamics  } 

\noindent 
In this section, we follow   the ideas proposed in 
\cite{RQ13,RVGSA13,VRFG13}.
The space of states $Q$ has a natural Riemannian structure. 
Therefore, it is natural to choose a cost function 
that is intrinsic and invariant, and in consequence non sensible to the change 
of coordinates. 
Following  Rojas Qinteros's thesis \cite{RQ13}, 
 we introduce a particular cost function  
to control the dynamics (\ref{control-evolution}):
\moneq \label{cout-hamilton} 
J(u) \,=\, {1\over2} \, \int_0^T M_{k \ell} (q) \, u^k \, u^\ell \, {\rm d}t \, . 
\monend
The controlled system (\ref{control-evolution}) (\ref{cout-hamilton})
is of type (\ref{edo}) (\ref{fonction-cout}) with 
\moneq \label{traduction}  
 \left\{ \begin{array}{rll}  \displaystyle \,\,  
 Y   & =\, \, \{ q^j \,, \, {\dot q}^j \} \,, \quad   
    & f \, =\,  \{ Y_2^j \,, \, 
- \Gamma^j_{k \ell} \,  {\dot q}^k \,{\dot q}^\ell - M^{j \ell} \, \partial_\ell V + u^j \} 
\,, \\ \displaystyle  \,\,  
 \lambda  & = \,  \{ u^k \} \,, \quad  
 & g \, =\,   {1\over2} \, M_{k \ell} (Y_1) \, u^k \, u^\ell \, . 
\end{array}  \right. \monend
%
%
The Pontryagin method introduces Lagrangre multipliers (or adjoint states)
$ \, p_j \,$ and $ \, \xi_j \,$  to form the Hamiltonian 
$ \, {\cal H}(Y, \, P, \, \lambda) \,$ function of state $ \, Y \,$ defined in 
(\ref{traduction}) and adjoint  $ \, P \, $ obtained by combining the two adjoint states:  
\moneq \label{etat-adjoint}  
P \, = \,  \{ p_j \,, \, \xi_j \} \,
 \monend
and $ \, \lambda \, = \,  \{ u^k \} \,$ as proposed in (\ref{traduction}). 
Taking into account (\ref{hamiltonien}), (\ref{traduction}) and (\ref{etat-adjoint}), we have:
\moneq \label{global-hamiltonien}  
 {\cal H}(Y, \, P, \, \lambda) \,=\, p_j \, {\dot q}^j 
+ \xi_j \, \big[ - \Gamma^j_{k \ell} \,  {\dot q}^k \,{\dot q}^\ell 
- M^{j \ell} \, \partial_\ell V + u^j \big] -  {1\over2} \, M_{k \ell} (Y_1) \, u^k \, u^\ell \, . 
 \monend

\bigskip \monitem   {\bf Proposition 4. \quad  Interpretation of one adjoint state.  }

\noindent 
When the cost function $J$ defined in (\ref{cout-hamilton}) is stationary, 
the adjoint state $ \, \xi^j \,$
is exactly equal to the applied  force (and torque!)  $ \, u^j \, $ in the right hand side of
the dynamic equation (\ref{control-evolution}): 
\moneq \label{interpretation}  
 \xi^j \,= \,  u^j \, . 
 \monend
%

%
\smallskip \noindent  {\bf Proof of Proposition 4.}

\noindent 
Due to the expression (\ref{global-hamiltonien}) of the Hamiltonian function, 
the optimality condition $ \, {{\partial {\cal H}}\over{\partial \lambda}} = 0 \,$ takes
the simple form 
\moneqstar
\xi_j \,= \, M_{j \ell} \, \xi^\ell  \, . 
\monendstar 
This relation is equivalent 
to the condition (\ref{interpretation}). \hfill $\square$ 

\smallskip \monitem
The reduced Hamiltonian $ \, H(Y,\, P) \,$ at the optimum can be 
explicited without difficulty. We just replace the control force $ \, u_j \,$
by the adjoint state $ \, \xi_j$: 
\moneqstar 
H(Y,\, P) \, = \, p_j \, {\dot q}^j 
+ \xi_j \, \big[ - \Gamma^j_{k \ell} \,  {\dot q}^k \,{\dot q}^\ell 
- M^{j \ell} (Y_1) \, \partial_\ell V \big] +  {1\over2} \,  M^{k \ell} (Y_1) \, \xi_k \, \xi_\ell \, . 
 \monendstar 
The symplectic dynamics (\ref{symplectic}) can be written simply:
\moneq \label{dyn-symplectic}  
{\dot q}^j \,=\, {{\partial H}\over{\partial p_j}} \,, \quad 
{\ddot q}^j \,=\, {{\partial H}\over{\partial \xi_j}} \,, \quad 
{\dot p}_j \,=\, -{{\partial H}\over{\partial q^j}} \,, \quad 
{\dot \xi}_j \,=\, -{{\partial H}\over{\partial {\dot q}^j}} \, . 
 \monend
The two first equations of (\ref{dyn-symplectic}) give the initial controlled
dynamics (\ref{control-evolution}).  We have also 
\moneqstar 
 \left\{  \begin{array}{l}  \displaystyle  \quad 
{{\partial H}\over{\partial q^j}} \, = \, 
- \big ( \partial_j \Gamma^i_{k \ell} \big) \,  {\dot q}^k \,{\dot q}^\ell \, \xi_i
- \partial_j \big( M^{i \ell} \, \partial_\ell V \big) \, \xi_i  
+  {1\over2} \big ( \partial_j M^{k \ell} \big) \, \xi_k \, \xi_\ell 
 \\ \displaystyle  \qquad 
{{\partial H}\over{\partial {\dot q}^j}} \, = \, p_j - 
2 \, \Gamma^i_{k j} \,  {\dot q}^k \, \xi_i \, . 
 \end{array} \right. \monendstar
%
%
%
We deduce the developed form of the two last equations of  (\ref{dyn-symplectic}):
\moneq \label{p-point-j}  
 {\dot p}_j  \, = \, \big( \partial_j \Gamma^i_{k \ell} \big) \,   {\dot q}^k \,{\dot q}^\ell \, \xi_i 
+ \partial_j \big( M^{i \ell} \, \partial_\ell V \big) \, \xi_i
- {1\over2} \, \big( \partial_j  M^{k \ell} \big) \, \xi_k \, \xi_\ell  
 \monend 
\moneq \label{xi-point-j}  
 {\dot \xi}_j  \, = \, 2 \, \Gamma^i_{k j} \,   {\dot q}^k \, \xi_i - p_j \, . 
 \monend

\bigskip  \noindent {\bf \large 5) \quad Intrinsic evolution of the generalized force } 

 \noindent
We introduce the covector $ \, \xi \,$ according to its covariant coordinates:
$ \,\,  \xi   \, = \,  \xi_j \,  \partial^j  $.
We have the following result, first established in \cite{RQ13} and \cite{VRFG13}:

\bigskip \monitem   {\bf Proposition 5. \quad  Covariant evolution equation of the optimal force.}

\noindent 
With the above notations and hypotheses, the forces and torques  $ \, u \,$ 
satisfy the following time evolution: 
\moneq \label{evolution-force}  
\Big( {{{\rm d}^2 u}\over{{\rm d}t^2}} \Big)_j +
R^i_{k \ell j} \,   {\dot q}^k \,{\dot q}^\ell \, u_i  + 
\big( \nabla^2_{j k} V \big) \, u^k \,=\, 0 \, . 
 \monend
%

%
\smallskip \noindent  {\bf Proof of Proposition 5.}

\noindent 
The time covariant derivative of the covector $ \, \xi \,$ is given by 
\moneqstar  
{{{\rm d} \xi}\over{{\rm d}t}} \,=\,  
\big( {\dot \xi}_j - \Gamma^i_{j k} \,   {\dot q}^k \,  \xi_i \big)  \, \partial^j
\monendstar
%
%
that is 
$ \displaystyle \,\, \Big(  {{{\rm d} \xi}\over{{\rm d}t}} \Big)_j \,=\, 
 {\dot \xi}_j - \Gamma^i_{j k} \,   {\dot q}^k \,  \xi_i  . \, $ 
 We report this expression in (\ref{xi-point-j}):

\moneq \label{p-j}  
 {p}_j  \, = \, \Gamma^i_{j k} \,  {\dot q}^k \,  \xi_i  
- \Big(  {{{\rm d} \xi}\over{{\rm d}t}} \Big)_j \, . 
 \monend
%
We wish to differentiate relative to time the expression $ \, p_j \, $  given in 
(\ref{p-j}).  The covariant derivatives of the covector $p$ can be evaluated as follows:
\moneqstar  
\Big( {{{\rm d} p}\over{{\rm d}t}} \Big)_j \,=\,  {\dot p}_j  - 
 \Gamma^\ell_{j k} \,  {\dot q}^k \, p_\ell  \, . 
\monendstar
Then we have, taking into account again the relation (\ref{p-j}):

\smallskip \noindent $ \displaystyle  
{\dot p}_j \,=\, {{{\rm d}}\over{{\rm d}t}}  \Big[ 
 \Gamma^i_{j k} \,  {\dot q}^k \,  \xi_i - \Big( {{{\rm d} \xi}\over{{\rm d}t}} \Big)_j \Big] 
+ \Gamma^\ell_{j k} \,  {\dot q}^k \, \Big[  \Gamma^i_{s \ell} \,  {\dot q}^s  \,  \xi_i 
-  \Big( {{{\rm d} \xi}\over{{\rm d}t}} \Big)_\ell \Big]  $ 

\smallskip \noindent $ \displaystyle \quad \,\,  = \, 
 \partial_\ell \big(  \Gamma^i_{j k} \big) \,  {\dot q}^k \,  {\dot q}^\ell \,  \xi_i 
+  \Gamma^i_{j k} \,  {\ddot q}^k \,  \xi_i  
+  \Gamma^i_{j k} \,  {\dot q}^k \,  \Big(  {{{\rm d} \xi}\over{{\rm d}t}} \Big)_i 
- \Big( {{{\rm d}^2 \xi}\over{{\rm d}t^2}} \Big)_j $

 \qquad \qquad  $ \displaystyle 
+  \Gamma^\ell_{j k}  \, \Gamma^i_{s \ell} \,  {\dot q}^k \,   {\dot q}^s \,  \xi_i
- \Gamma^\ell_{j k} \,  {\dot q}^k \, \Big( {{{\rm d} \xi}\over{{\rm d}t}} \Big)_\ell  $


\smallskip \noindent $ \displaystyle \quad = \, 
 \partial_\ell \big( \Gamma^i_{j k} \big) \,   {\dot q}^k  \,  {\dot q}^\ell  \, \xi_i 
+ \Gamma^s_{k j} \, \Gamma^i_{s \ell} \,   {\dot q}^k  \,  {\dot q}^\ell  \, \xi_i 
- \Big( {{{\rm d}^2 \xi}\over{{\rm d}t^2}} \Big)_j 
+ \Gamma^i_{k j} \, {\ddot q}^k \, \xi_i   $ 

\hfill 
due to the simplification of two terms
%
%


\smallskip \noindent $ \displaystyle \quad = \, 
 \partial_\ell \big( \Gamma^i_{j k} \big) \,   {\dot q}^k  \,  {\dot q}^\ell  \, \xi_i 
+ \Gamma^s_{k j} \, \Gamma^i_{s \ell} \,   {\dot q}^k  \,  {\dot q}^\ell  \, \xi_i 
- \Big( {{{\rm d}^2 \xi}\over{{\rm d}t^2}} \Big)_j 
$

\qquad \qquad  $ \displaystyle 
 + \Gamma^i_{k j} \, \Big( - \Gamma^k_{\ell s} \,  {\dot q}^s  \,  {\dot q}^\ell
- M^{k \ell} \, \partial_\ell V  + \xi^k \Big) \, \xi_i $
\qquad \qquad  due to (\ref{control-evolution})

\smallskip \noindent $ \displaystyle \quad = \, 
\Big(   \partial_\ell  \Gamma^i_{j k} + \Gamma^s_{j k} \, \Gamma^i_{s \ell} 
-  \Gamma^s_{k \ell} \, \Gamma^i_{s j} \Big) \,  {\dot q}^k  \,  {\dot q}^\ell  \, \xi_i 
- \Big( {{{\rm d}^2 \xi}\over{{\rm d}t^2}} \Big)_j 
- \Gamma^i_{k j} \, M^{k \ell} \, \partial_\ell V \, \xi_i 
+ \Gamma^i_{k j} \, \xi^k \, \xi_i  $

 \noindent and 
\moneq \label{p-point-j-seconde}  
 {\dot p}_j  \, = \,
\Big( R^i_{k j \ell} + \partial_j \Gamma^i_{k \ell}  \Big) \,  {\dot q}^k  \,  {\dot q}^\ell  \, \xi_i 
- \Big( {{{\rm d}^2 \xi}\over{{\rm d}t^2}} \Big)_j 
- \Gamma^i_{k j} \, M^{k \ell} \,  \partial_\ell V \, \xi_i 
+ \Gamma^i_{k j} \, \xi^k \, \xi_i  
\monend 
taking into account the expression (\ref{tenseur-riemann})  of the Riemann tensor.
We confront the relations (\ref{p-point-j-seconde}) and  (\ref{p-point-j}). We
deduce 
\moneq \label{intermediaire}  
 \left\{  \begin{array}{l}  \displaystyle  
R^i_{k j \ell} \,  {\dot q}^k  \,  {\dot q}^\ell  \, \xi_i 
- \Big( {{{\rm d}^2 \xi}\over{{\rm d}t^2}} \Big)_j 
- \Gamma^i_{k j} \, M^{k \ell} \, \partial_\ell V \, \xi_i 
+ \Gamma^i_{k j} \, \xi^k \, \xi_i  \\ \displaystyle 
\qquad \qquad \qquad = \,\, \partial_j \big( M^{i \ell} \, \partial_\ell V \big) \,  \xi_i 
- {1\over2} \, \big( \partial_j  M^{k \ell} \big) \, \xi_k \, \xi_\ell \, . 
 \end{array} \right. \monend
%
%
We take into account the second Ricci identity (\ref{ricci}). It comes

\smallskip \noindent $$ \displaystyle 
 \big( \partial_j  M^{k \ell} \big) \, \xi_k \, \xi_\ell \,=\, 
- \Gamma^k_{j s} \,  \xi_k \, \xi^s  - \Gamma^\ell_{j s} \,  \xi_\ell \, \xi^s
\,=\, -2 \, \Gamma^k_{j \ell} \,  \xi_k \, \xi^\ell \,=\,  
-2 \, \Gamma^i_{j k} \,  \xi^k \, \xi_i   \, . $$ 

\smallskip \noindent 
Then we can write the relation (\ref{intermediaire}) in a simpler way:

\smallskip \noindent $ \displaystyle 
R^i_{k j \ell} \,  {\dot q}^k  \,  {\dot q}^\ell  \, \xi_i 
  - \Big( {{{\rm d}^2 \xi}\over{{\rm d}t^2}} \Big)_j \,=\, 
\Big[ \Gamma^i_{k j} \, M^{k \ell} \, \partial_\ell V
+ \partial_j \Big( M^{i \ell} \, \partial_\ell V \Big) \Big] \, \xi_i  $
 
\smallskip \noindent $ \displaystyle \qquad = \, 
\Big[ \Gamma^i_{k j} \, M^{k \ell} \, \partial_\ell V
- \Gamma^i_{j s} \, M^{s \ell} \, \partial_\ell V 
- \Gamma^\ell_{j s} \, M^{i s} \, \partial_\ell V 
+ M^{i \ell} \,  \partial_\ell \partial_j V \Big] \, \xi_i  $

\smallskip \noindent $ \displaystyle \qquad = \, 
\Big(  \partial_\ell \partial_j V - \Gamma^s_{j \ell} \,  \partial_s V \Big) \,  M^{i \ell} \, \xi_i  
\,\, = \,\,  \big(  \nabla^2_{j \ell} V  \big) \, \xi^\ell  $

\smallskip \noindent 
due to the expression (\ref{second-gradient})  of the second gradient of a scalar field. 
We have established the following evolution equation
\moneqstar
\Big( {{{\rm d}^2 \xi}\over{{\rm d}t^2}} \Big)_j +
\big( \nabla^2_{j k} V \big) \, \xi^k \,=\, 
R^i_{k j \ell} \,   {\dot q}^k \,{\dot q}^\ell \, \xi_i 
 \monendstar
and the relation (\ref{evolution-force}) is a simple consequence 
of the anti-symmetry of the Riemann tensor and of the identity 
(\ref{interpretation}). 
\hfill $\square$

\bigskip   \noindent {\bf Conclusion }

 \noindent
We have established that the methods of Euler-Lagrange and Pontryagin  
conduct to two second order differential systems that couples state 
and control variables. 
The choice of a Riemannian metric allows the two systems to be in a well-defined 
tensorial nature: contravariant for the equation of motion and covariant 
for the equation of the control variables. 
The study of a robotic system, of which we try to optimize the control, 
shows how important is the introduction of an appropriate geometric structure. 
Riemannian geometry selected on the configuration parameter space favors the metric 
directly related to the mass tensor  as suggested by the expression of the kinetic energy.
An undeniable impact is the choice of an invariant cost function
with respect  to the choice of parameters, 
this is a stabilizing factor for numerical developments.
Pontryagin's principle applied to contravariant equation of motion associated 
with the cost function conducts to  a mechanical  interpretation of  adjoint states. 

 \noindent
The adjoint control equation  is established in a condensed form by the introduction 
of  second order covariant  derivatives  and shows the Riemann curvature tensor. 
Moreover, this framework exhibits a numerically stable method when discretization 
is considered.  The resolution of the coupled system gives a direct access 
to control variables without any additional calculation.
Thus, future numerical  developments will have to juggle between two coupled systems 
of second-order ordinary differential equations: the equation of motion 
and the equation for the control.

%




\begin{thebibliography}{}

\end{thebibliography}


\bigskip \bigskip   \noindent {\bf \large  References } 

 \vspace{-.2cm}  \noindent 
\begin{thebibliography}{99}




\bibitem{LV95} {\vskip -.4cm} 
M. Lazrak, C. Vall\'ee.
``Commande de robots en temps minimal'', 
{\it Revue d'Automatique et de Productique Appliqu\'ees 
(RAPA)}, volume~8, issue 2-3, p.~217-222, 1995. 


\bibitem{LR75} 
D. Lovelock, H. Rund.
{\it Tensors, Differential Forms and Variational Principles}, 
John Wiley \&  Sons, New York, 1975. 


\bibitem{PBGM62} 
L.S. Pontryagin, V.G. Boltyanskii, R.V.  Gamkrelidze, E.F. Mishchenko.
{\it The Mathematical Theory of Optimal Processes}, 
english translation. Interscience, 1962. 



\bibitem{RQ13} 
J.A. Rojas Quintero.
{\it Contribution \`a la manipulation dextre dynamique pour les aspects 
conceptuels et de commande en ligne optimale}, Thesis
Poitiers University, 31 October 2013. 


\bibitem{RVGSA13} 
J.A. Rojas Quintero, C. Vall\'ee, J.P. Gazeau, P. Seguin, M. Arsicault.
``An alternative to Pontryagin's principle 
for the optimal control of jointed arm robots'', 
Congr\`es Fran\c cais de M\'ecanique, Bordeaux, 26 - 30 August 2013. 

\bibitem{Si12} 
R. Siebert.
{\it  Mechanical integrators for the optimal control in multibody dynamics}, 
Dissertation, Department Maschinenbau. Universit\"at Siegen, 2012. 


\bibitem{VRFG13} 
C. Vall\'ee,  J.A. Rojas Quintero,  D. Fortun\'e,  J.P. Gazeau. 
``Covariant formulation of optimal control of jointed arm robots: 
an alternative to Pontryagin's principle'', 
arXiv:1305.6517, 28 May 2013.






\end{thebibliography}
\end{document}